\documentclass[12pt]{article}
\usepackage{amsmath}
\usepackage{graphicx,psfrag,epsf}
\usepackage{enumerate}
\usepackage{natbib}
\usepackage{amssymb}
\usepackage[colorlinks,citecolor=blue,urlcolor=blue]{hyperref}
\usepackage{mathrsfs}
\usepackage{amsfonts,anysize,amsmath,indentfirst,geometry,xcolor}
\usepackage{bm}
\usepackage{amsthm}
\usepackage{multirow}
\usepackage{chngpage}
\usepackage{array}
\usepackage{rotating}
\usepackage{setspace}
\usepackage{graphicx}
\usepackage{algorithm}
\usepackage{algorithmic}
\usepackage{epsfig}
\usepackage{amsbsy}

\newtheorem{remark}{Remark}

\newcommand{\blind}{0}

\begin{document}
\doublespacing

\def\spacingset#1{\renewcommand{\baselinestretch}%
{#1}\small\normalsize} \spacingset{1}

 \date{}
\if0\blind
{
  \title{\bf Jackknife Partially Linear Model Averaging for the Conditional Quantile Prediction}
  \author{Jing Lv\\
    School of Mathematics and Statistics,\\
     Southwest University, Chongqing, China}
  \maketitle
} \fi

\if1\blind
{
  \bigskip
  \bigskip
  \bigskip
  \begin{center}
    {\LARGE\bf Title}
\end{center}
  \medskip
} \fi

\bigskip

\begin{abstract}
Estimating the conditional quantile of the interested variable with respect to changes in the covariates is frequent in many economical applications as it can offer a comprehensive insight. In this paper, we propose a novel semiparametric model averaging to predict the conditional quantile even if all models under consideration are potentially misspecified. Specifically, we first build a series of non-nested partially linear sub-models, each with different nonlinear component. Then a leave-one-out cross-validation criterion is applied to choose the model weights. Under some regularity conditions, we have proved that the resulting model averaging estimator is asymptotically optimal in terms of minimizing the out-of-sample average quantile prediction error. Our modelling strategy not only effectively avoids the problem of specifying which a covariate should be nonlinear when one fits a partially linear model, but also results in a more accurate prediction than traditional model-based procedures because of the optimality of the selected weights by the cross-validation criterion. Simulation experiments and an illustrative application show that our proposed model averaging method is superior to other commonly used alternatives.
\end{abstract}

\noindent%
{\it Keywords:} Asymptotic optimality, B-splines, Conditional quantile prediction, Leave-one-out cross-validation, Model averaging, Partially linear models.

\spacingset{1.45}

\section{Introduction}\label{sec:1}

In many situations of practical interest, especially for econometrics, social sciences and medical fields,
we are more concerned to predict the conditional quantiles of interested variables because a full range of quantile analysis provides a broader insight than the classical mean regression \citep{K05}. For example, in business and economics, petroleum is a primary source of nonrenewable energy, and has important
influence on industrial production, electric power generation, and transportation. Most economists take care of the high quantiles of oil prices, because oil price fluctuations have considerable effects on economic activity. In the past decades, traditional parametric and semiparametric modelling strategies for quantile regression have been well developed, including \cite{k07,bc2011,KLZ11,wwl2012,fz2016,fb16,mh16,fbf21} and among others.

In practice, the underlying model is often unknown and all models under consideration are potentially misspecified. We all know that it is difficult to find a optimal model for a dataset of interest.
Model averaging, as a well-known ensemble technique, combines a set of candidate
models by assigning heavier weights to stronger models. The main superiority of model averaging is that it effectively incorporates useful information from all possible candidate models and thus substantially reduces the risk of misspecification and generally yields more accurate prediction results than a single selected model. For example, to explain a specific economic phenomenon, many plausible candidate models are all useful. In that case, using an averaged model instead of a particular model, the risk arising from misspecification can be reduced markedly. Over the past decade or so, model averaging for condition mean regression has been developed rapidly, see \cite{wzz10,hr12,AL14,zl18,zwlc20,fl21}. However, limited works have been done for studying quantile model averaging. Recently, \cite{LS15} introduced a jackknife quantile model averaging procedure with the optimal model weights by minimizing a leave-one-out cross-validation criterion. \cite{wz19a} introduced a jackknife model averaging for composite quantile regression, which can be regarded as an extension of \cite{LS15}. Instead of choosing the optimal model weights, \cite{ls21} proposed to average over the complete subsets for quantile regression, where the optimal size of the complete subset is selected by the cross-validation.

So far the mentioned above works mainly focus on averaging a set of parameterized models such as linear regression models. We have to acknowledge that such simple models are easy to interpret and widely accepted by scientific researchers. However, it is understood that, in practice, the response variable may depend on the predictors in a very complicated manner. If only parametric sub-models are adopted, it is hard for us to acquire satisfactory prediction results because they fail to capture the complicated relationship between the response and predictors. Though all candidate models might be misspecified in reality, we hope that the approximation
capabilities might be improved by using more flexible semiparametric sub-models.
An alternative approach, one considers here, is to construct a weighted average of a series of flexible semiparametric sub-models. We may refer to \cite{llrz18,lxwn18,zl18,zwzz19,zw19} for reviews of recent developments on semiparametric model averaging. However, these aforementioned research findings only concern the conditional mean prediction, and discussions on quantile regression are rather limited.

The partially linear model (PLM) introduced by \cite{egrw86}, as one of the most popular semiparametric models, has been received extensive attention due to its flexible specification.
The main merit of the PLM is that it does not require the parametric assumption for all covariates and allows one to capture potential nonlinear
effects. Although we have witnessed a booming development of the PLM in recent years (e.g., \cite{h00,lwc07,ll09,xh09,zcl11}), these methodologies are all based on the assumption that a correctly specified model is given. So far little work has been done on quantile model averaging for the PLM.
In this work, we will exploit a semiparametric model averaging by optimally combining a series of PLMs to achieve the goal of flexible conditional quantile prediction. This fills an important gap in semiparametric model averaging for the conditional quantile prediction.

The contribution of this paper is three folds. First, it is usually a challenging job to decide which a covariate should be nonlinear when one fits a PLM. Actually, any continuous covariate can be taken as the nonparametric component. Our proposed model averaging effectively avoids the criticism of
artificially specifying the nonparametric component in PLMs because we average multiple partially linear sub-models (each with different nonlinear component) by assigning heavier weights to stronger sub-models.
An another superiority of the proposed approach is that it is more robust against model misspecification
, and thus outperforms than traditional model-based approaches (e.g., linear models, partially linear models and additive models) and parametric model averaging procedures.
Second, we rapidly estimate the entire conditional quantile process over $(0,1)$ of the interested response rather than a discrete set of quantiles by modeling quantile regression coefficients as parametric functions of quantile level. Compared with standard quantile estimation procedure, the strategy of modeling quantile functions parametrically simplifies calculation and gains better estimation efficiency because of utilizing all useful information across quantiles.
Third, we prove that the proposed model averaging estimator is asymptotically optimal in the sense that its out-of-sample average quantile prediction error is asymptotically identical to that of the best but infeasible model averaging estimator. It is instructive to mention that our theoretical results intrinsically distinguish from those of \cite{LS15},\cite{wz19a} and \cite{ls21} who focus on parametric model averaging.

\section{Methodology}\label{sect.2}
\subsection{Model and Estimation}\label{sect.2.1}

Let $\left\{ \left(Y_i,\bm{X}_i\right) \right\}_{i=1}^n$ be independent and identically distributed samples of $\left(Y,\bm X\right)$ with $n$ individuals, where $Y$ is a scalar response variable and $\bm{X}$ is the vector of covariates. For a given quantile level $\tau\in (0, 1)$, let ${Q}\left({Y}| \bm X,\tau   \right)$ be the $\tau$th conditional quantile function of the response $Y$ given the covariates $\bm{X}$. Without loss of generality, the covariates $\bm{X}$ are allowed to be discrete or continuous. Suppose that $\bm{X}=\left(\bm X_{\cal A}^\top,\bm X_{\cal B}^\top\right)^\top$, where $ \bm X_{\cal A}$ and $\bm X_{\cal B}$ are vectors of $p$-dimensional continuous and $q$-dimensional discrete variables respectively, and $^\top$ is the transpose of a vector or matrix. Our goal is to estimate ${Q}\left({Y}| \bm X,\tau   \right)  \triangleq \mu\left(\bm X, \tau\right)$, which is of particular use for prediction. This
is also the typical goal in the optimal model averaging literature (\cite{LS15,zw19,ls21}).

As far as we know, it is infeasible to estimate $\mu\left(\bm X, \tau\right)$ without any structure assumption for $p+q>2$ because of the curse of dimensionality. Partially linear models \citep{h00,lwc07,ll09}, as one of the most commonly used semiparametric models, have been
developed to resolve this problem due to their flexible specification. However, all models under investigation might be incorrect in practice. Using a single model might ignore the useful information from the other models, and thus results in a poor predictive performance. As an attractive alternative, aggregating candidate models with a weighted average effectively provides a better approximation to the true quantile function $\mu\left(\bm X, \tau\right)$ and gives us great potentials for future prediction.

Specifically, we assume that each candidate model has the partially linear model structure.
In practice, one will encounter the uncertainty of whether a covariate should be in the linear or nonlinear given that it is in the model. In fact, any continuous elements of $\bm X$ might be taken as the nonparametric components. To avoid the criticism of artificially deciding which covariates are nonlinear in PLMs, we build a sequence of $p$ partially linear sub-models $\mathbb{M}_1,\cdots,\mathbb{M}_p$, and the model weights automatically adjust the relative importance of these sub-models. By taking the $s$th element of $\bm X_i$, $X_{is}$ as the nonparametric component, the $s$th sub-model $\mathbb{M}_s$ is given by
\begin{align} \label{eq1}
\mathbb{M}_s:
 {\mu ^{(s)}}\left(\bm X_i,\tau\right) \triangleq  g^{(s)}\left( X_{is},\tau\right)+\bm X_{i\setminus s}^\top\bm \beta^{(s)}(\tau)  ,  i=1,\cdots,n,s=1,\cdots,p,
\end{align}
where $\bm X_{i\setminus s}$ is a $(p+q-1)\times 1$ covariate vector by removing the $s$th predictor $X_{is}$,
$\bm\beta^{(s)}(\tau)=\left(\beta_j^{(s)}(\tau): 1\leq j \leq p+q-1 \right)^\top$ and $ g^{(s)}\left( \cdot, \tau\right)$ are unknown parameter vector and smooth function at the $\tau$-th quantile. Note that ${\mu ^{(s)}}\left(\bm X_i,\tau\right) $ is the condition quantile function under the $s$th sub-model $\mathbb{M}_s$. Although the ``intercept'' term does not appear in model (\ref{eq1}), it is actually included in the functional component.
It is easy to see that the differences between any two candidate models lie, not only in linear components, but also in choosing which a covariate should be taken as the nonparametric element. To offer an optimal weighting scheme, we first should estimate the unknown parameter vector $\bm \beta^{(s)}(\tau)$ and function $g^{(s)}\left( \cdot,\tau\right)$ of each candidate model.

To estimate the functional component $g^{(s)}\left( \cdot,\tau\right)$, we can approximate $g^{(s)}\left( \cdot,\tau\right)$ by B-spline basis functions because of its efficient in function approximation and stable numerical computation (see \cite{De01}). Under proper conditions on $g^{(s)}\left( \cdot,\tau\right)$ (e.g., Condition (C2) below), according
to Corollary 4.10 of \cite{S81}, we can approximate $g^{(s)}\left( \cdot,\tau\right)$ as
\begin{align} \label{eq2}
g^{(s)}\left( \cdot,\tau\right) \approx\bm B^\top\left(\cdot\right) \bm\gamma^{(s)}(\tau),
\end{align}
where $\bm B\left(\cdot\right) = \left( B_c\left(\cdot\right),:1 \leq c \leq J_n\right)^\top$ is a $J_n \times 1$ vector of normalized B-spline basis functions of order $d$ ($d \geq 2$), $J_n=N_n+d$, $N_n$ is the number of interior knots and $\bm\gamma^{(s)}(\tau)= \left(\gamma_j^{(s)}(\tau): 1 \leq j \leq J_n \right)^\top$. Equally spaced knots are used here for technical simplicity. However other regular
knot sequences can also be used, with similar asymptotic results.
Then, substituting (\ref{eq2}) into the model (\ref{eq1}), we can get
\begin{align} \label{eq3}
\mathbb{M}_s: {\mu^{(s)}}\left(\bm X_i,\tau\right)
 \approx &\bm B^\top\left(X_{is}\right) \bm\gamma^{(s)}(\tau) +\bm X_{i\setminus s}^\top\bm \beta^{(s)}(\tau) \nonumber\\
\triangleq &  \bm Z_i^{(s)\top} \bm\xi^{(s)}(\tau),
\end{align}
where $\bm Z_i^{(s)}=\left(\bm B^\top\left(X_{is}\right), \bm X_{i\setminus s}^\top\right)^\top$ and $\bm\xi^{(s)}(\tau)=\left(\bm\gamma^{(s)\top}(\tau), \bm \beta^{(s)\top}(\tau) \right)^\top$.

Obviously, $\bm\xi^{(s)}(\tau)$ can be regarded as a set of quantile regression coefficient functions
describing how each regression coefficient depends on the quantile level $\tau$. We might obtain an estimator of $\bm\xi^{(s)}(\tau)$ at a single quantile of interest by using the standard quantile regression (e.g., \cite{K05}). However, to improve efficiency of the coefficients estimates, we adopt the strategy of \cite{fb16} to model $\bm\xi^{(s)}(\tau)$ by a series of parametric functions.
Specifically, we take $\bm\xi^{(s)}(\tau)$ as a function of quantile level $\tau$ that relies on a finite-dimensional parameter $\bm\theta^{(s)}$
\begin{align} \label{eq4}
\bm\xi^{(s)}(\tau)=\bm\theta^{(s)} \bm b\left(\tau\right),
\end{align}
where $\bm b\left(\tau\right)=\left(b_j\left(\tau\right): 1\leq j \leq K \right)^\top$ is a set of $K$ known basis functions of $\tau \in (0,1)$, and $\bm\theta^{(s)}$ is a $(J_n+p+q-1) \times K$ matrix with the $(u,v)$th entries $\theta_{uv}^{(s)}$ for $u=1,\cdots,J_n+p+q-1 $ and $v=1,\cdots,K $. Under the model \eqref{eq4}, we have $\xi _l^{(s)}(\tau ) = \sum\limits_{k = 1}^K {{\theta _{lk} ^{(s)}}} {b_k}\left( \tau  \right),l = 1,\cdots ,J_n+p+q-1$. In practice, to obtain an estimate of $\bm\theta^{(s)}$, we need to specify $\bm b\left(\tau\right)$ in advance.
As mentioned in \cite{fb16,YCC17,fbf21}, valid choices of $\bm b\left(\tau\right)$ are, for example, functions of the form $\tau^{\alpha}$, $\log(\tau)$, $\log(1-\tau)$, $\alpha^\tau$, the quantile function of any distribution
with finite moments, splines, or a combination of the above. In general, the selected basis set $\bm b\left(\tau\right)$ should satisfy two conditions. First, $\bm Z_i^{(s)\top} \bm\theta^{(s)} \bm b\left(\tau\right)$ defines a valid quantile function (i.e., is an increasing function of $\tau$) for some $\bm\theta^{(s)}$ at the observed values of $\bm Z_i^{(s)}$. Second, it is differentiable in the interior of its support. The simulation results of Table \ref{table1} show that the proposed method is not sensitive to the selection of $\bm b\left(\tau\right)$.

To facilitate the presentation, we need to introduce some notations.
Let ${\rm{Vec}}\left( \cdot \right)$ be the vectoring operation, which creates a column vector by stacking the column vectors of below one another, that is, ${\rm{Vec}}\left( \bm\theta^{(s)} \right)=(\bm \theta_1^{(s)\top},\cdots,\bm \theta_K^{(s)\top})^\top$, where $\bm \theta_k^{(s)}=\left(\theta_{1,k}^{(s)},\cdots,\theta_{J_n+p+q-1,k}^{(s)}\right)^\top$ is the $k$-th column of the parameter matrix $\bm \theta^{(s)}$ in the $s$-th candidate model. Define $\bm D_i^{(s)}\left( \tau  \right) = \bm b\left( \tau  \right) \otimes \bm Z_i^{(s)}$ with $\otimes$ representing the Kronecker product of two matrices. Then, we have
 $ \bm Z_i^{(s)\top}\bm\theta^{(s)} \bm b\left(\tau\right) =\bm D_i^{(s)\top}\left( \tau  \right){\rm{Vec}}\left( \bm\theta^{(s)} \right)$ and $\bm\xi^{(s)}(\tau)= \sum\limits_{k = 1}^K {\bm \theta_k^{(s)}b_k\left(\tau\right) } $.

Motivated by \cite{fb16}, we can further integrate information from different quantile levels to improve efficiency and obtain the estimator ${\rm{Vec}}\left( \hat{\bm\theta}^{(s)} \right)$ of ${\rm{Vec}}\left( \bm\theta^{(s)} \right)$ by minimizing the integrated loss function
\begin{align}\label{eq5}
{\rm{Vec}}\left( \hat{\bm\theta}^{(s)} \right)=\mathop {\arg \min }\limits_{{\rm{Vec}}\left( \bm\theta^{(s)} \right)} \bar{\cal L}_n\left( {\rm{Vec}}\left( \bm\theta^{(s)} \right) \right)=\mathop {\arg \min }\limits_{{\rm{Vec}}\left( \bm\theta^{(s)} \right)}\int_0^1 {\cal L}_n\left( {\rm{Vec}}\left( \bm\theta^{(s)} \right) \right) d\tau,
\end{align}
where ${\cal L}_n\left( {\rm{Vec}}\left( \bm\theta^{(s)} \right)\right)=n^{-1} \sum\limits_{i = 1}^n {\rho_{\tau} \left( Y_i-\bm D_i^{(s)\top}\left( \tau  \right){\rm{Vec}}\left( \bm\theta^{(s)} \right)\right)}$ and ${\rho _\tau }\left( u \right) = u\left( {\tau  - I\left( {u < 0} \right)} \right)$ is the quantile check function. The objective function $\bar{\cal L}_n\left( {\rm{Vec}}\left( \bm\theta^{(s)} \right) \right)$ can be regarded as an average loss function, achieved by marginalizing ${\cal L}_n\left( {\rm{Vec}}\left( \bm\theta^{(s)} \right) \right)$ over the  entire interval $(0, 1)$. In addition, the solution of minimizing \eqref{eq5} is currently implemented by the {\tt iqr} function in the {\tt qrcm} R package. Define $\bm{\hat\xi}^{(s)}(\tau)=\sum\limits_{k = 1}^K {\hat{\bm \theta}_k^{(s)}b_k\left(\tau\right) } $ and ${\hat\mu ^{(s)}}\left(\bm X_i,\tau\right)=\bm D_i^{(s)\top}\left( \tau  \right){\rm{Vec}}\left( \hat{\bm\theta}^{(s)} \right)$ for $s=1,\cdots,p$.

With the estimators of each candidate model readily available, the model averaging estimator of $\mu\left(\bm X_i, \tau\right)$ is thus expressed as
\begin{align}\label{eq6}
{\hat \mu }^{\bm w}\left( {{\bm X_i},\tau } \right) = \sum\limits_{s = 1}^p {{{ w}_s}} \hat \mu ^{(s)}\left( {{\bm X_i},\tau } \right),
\end{align}
where $\bm w=(w_j,1\leq j \leq p)^\top$ is the model weight vector belonging to the set
\[\mathbb{W}=\left\{ \bm w\in [0,1]^{p}:  \sum\limits_{s=1}^{p} {w_s}=1 \right\} .\]

\begin{remark}\label{remark0}

The main merits of parametric modeling of $\bm\xi^{(s)}(\tau)$ include the following two aspects. One the one hand, the model \eqref{eq4} extracts the common features of the quantile regression coefficients $\bm\xi^{(s)}(\tau)$ over $\tau\in (0,1)$ via the $K$-dimensional known basis function vector $\bm b\left(\tau\right)$. Moreover, it permits estimating the entire quantile process rather than only obtaining a discrete set of quantiles. Thus, this modelling strategy presents numerous superiorities including a simpler computation, increased statistical efficiency and easy interpretability of the results. On the other hand, our proposed modelling strategy can effectively estimate bivariate functions $g^{(s)}\left( \cdot,\tau\right)$ by a combination of B-spline approximation and parametric modeling of $\bm\xi^{(s)}(\tau)$.
\end{remark}

\begin{remark}\label{remark1}
In practice, researchers are ignorant of the true model. All considered candidate models might be wrong, but each candidate model may characterize only some of the properties of the true data generating process. Although our constructed partially linear sub-models are more sophisticated and flexible than traditional linear models, they still may not be the true model. To reduce the risk of model misspecification, we construct a model average strategy to achieve accurate prediction for the conditional quantile function $\mu\left(\cdot, \tau\right)$ by assigning higher weights to the better sub-models. Compared with existing strategies of constructing nested sub-models (i.e., \cite{wzz10,hr12,LS15,zl18,zw19,zwlc20}), it is worth noting that each candidate model that we consider includes all covariates and thus our modelling strategy for sub-models is non-nested, which may be another attractive scheme of building semiparametric sub-models.
\end{remark}

\subsection{Jackknife Weighting}\label{sect.2.2}

Actually, the weight vector $\bm w$ in ${\hat \mu }^{\bm w}\left( {{\bm X_i},\tau } \right)$ is usually unknown and should be properly estimated as the choice of the weight plays a central role for model averaging
strategy. Following the idea of \cite{hr12,LS15}, we will adopt jackknife selection of $\bm w$ (also known as leave-one-out cross-validation). More specifically, we measure the average prediction error by the integrated loss minimization \citep{fb16} and define the leave-one-out cross-validation criterion by
\begin{equation} \label{eq7}
{\rm{CV}}_n\left( \bm  w \right) =\int_0^1 n^{-1} \sum\limits_{i = 1}^n {\rho_{\tau} \left( Y_i- \sum\limits_{s = 1}^p{w_s \bm D_i^{(s)\top}\left( \tau  \right){\rm{Vec}}\left( \hat{\bm\theta}^{(s)}_{[-i]} \right) } \right)} d\tau,
\end{equation}
where ${\rm{Vec}}\left( \hat{\bm\theta}^{(s)}_{[-i]} \right)$ is the jackknife estimator of ${\rm{Vec}}\left( \bm\theta^{(s)} \right)$ for $s=1,\cdots,p$, which is obtained by \eqref{eq5} without using the $i$th sample.
Minimizing ${\rm{CV}}_n\left( \bm  w \right)$ with respect to $\bm w$ leads to
\begin{equation} \label{eq8}
\hat{\bm w}= \mathop {\arg \min }\limits_{\bm w\in \mathbb{W}} {\rm{CV}}_n\left( \bm  w \right).
\end{equation}
Substituting $\hat {\bm w}$ for $\bm w$ in \eqref{eq6} results in the proposed model averaging estimator
\begin{align}\label{eq9}
{\hat \mu }^{\hat{\bm w}}\left( {{\bm X_i},\tau } \right) = \sum\limits_{s = 1}^p {\hat{w}_s} \hat \mu ^{(s)}\left( {{\bm X_i},\tau } \right)=\sum\limits_{s = 1}^p{\hat w_s \bm D_i^{(s)\top}\left( \tau  \right){\rm{Vec}}\left( \hat{\bm\theta}^{(s)} \right)}.
\end{align}
Averaging using the weight choice is called the \emph{jackknife quantile partially linear model averaging} ({\sf JQPLMA}).

Notice that the proposed estimator $\hat{\bm w}$ is computationally challenging due to
the complicated integrated loss function and constraint conditions on the weight vector.
Invoking the precursor work of \cite{KX14} and \cite{YCC17}, for the sake of computational convenience, we can approximate the objective function (\ref{eq7}) by
\begin{equation} \label{eq***}
{\rm{CV}}_n\left( \bm  w \right) \approx n^{-1}\sum\limits_{k = 1}^n n^{-1} \sum\limits_{i = 1}^n {\rho_{\tau_k} \left( Y_i- \sum\limits_{s = 1}^p{w_s \bm D_i^{(s)\top}\left( \tau_k  \right){\rm{Vec}}\left( \hat{\bm\theta}^{(s)}_{[-i]} \right) } \right)},\nonumber
\end{equation}
where $\tau_k=k/(n+1),k=1,\cdots,n$. Then a well-known nonlinear optimization such as the augmented
Lagrange method is considered to solve the constrained optimization problem, which is easy to be implemented by many software packages (i.e., the {\tt Rsolnp} in R).

Let $\left(\mathrm{y},\bm{\mathrm{x}}\right)$ be an independent copy of $(Y,\bm X)$. Write ${\mathcal{D}_n} = \left\{ {\left( {{Y_i},{\bm X_i}} \right)} \right\}_{i = 1}^n$, $\bm{\mathrm{x}}_{\setminus s}=\left(\mathrm{x}_j: 1 \leq j\leq p,j\neq s\right)^\top$, $\bm{\mathrm{z}}^{(s)}=\left(\bm B^\top\left(\mathrm{x}_s\right), \bm{\mathrm{x}}_{\setminus s}^\top\right)^\top$, $\mathbb{D}^{(s)}\left( \tau  \right) = \bm b\left( \tau  \right) \otimes \bm{\mathrm{z}}^{(s)}$ and $\hat \mu ^{(s)}\left( \bm{\mathrm{x}},\tau\right)= \mathbb{D}^{(s)\top}\left( \tau  \right){\rm{Vec}}\left( \hat{\bm\theta}^{(s)} \right)$ for $s=1,\cdots, p$. Define the out-of-sample average quantile prediction error (denoted as ${\rm OAQPE}_n$) as follows
\begin{align} \label{eq10}
{\rm{OAQPE}}_n\left( \bm w \right)=&E\left\{ \int_0^1{{\rho _\tau }\left(\mathrm{y} - \sum\limits_{s=1}^{p} {w_s\mathbb{D}^{(s)\top}\left( \tau  \right){\rm{Vec}}\left( \hat{\bm\theta}^{(s)} \right)} \right)d\tau\left| {{\mathcal{D}_n}} \right.} \right\}.
\end{align}
Next we will show that the weight vector selected by (\ref{eq8}) is asymptotically optimal in the sense of achieving the lowest possible ${\rm{OAQPE}}_n\left( \bm w \right)$ under some regularity conditions.

\begin{remark}\label{remark3}
What should be pointed out here is that the good in-sample performance does not necessarily indicate good out-of-sample performance because the future prediction of underlying models is partially or completely unknown to the practical users. Thus, we estimate the optimal weight vector by minimizing ${\rm{CV}}_n\left( \bm  w \right)$ instead of directly using the integrated loss function to obtain good out-of-sample prediction performance. In addition, it is understood that the selection of the loss function is closely related to the characteristic of the response variable's distribution that one wants to predict.
For example, the traditional quadratic (or quantile) loss function corresponds to the conditional mean (or quantile) of the distribution of the response. Here the object of our interest is the average of quantile prediction over the interval $(0,1)$, and thus it is natural to define the risk function (\ref{eq10}) which can be regarded as a beneficial extension of the criterion (2.13) in \cite{LS15}.
\end{remark}

\section{ Simulation Studies}\label{sect.4}

In this section, we conduct Monte Carlo experiments to examine the performance of
the proposed model average prediction procedure. To make a full comparison, we compare our proposal with the following popular model-based and model averaging prediction methods.

{\sf QLRM}: The traditional quantile linear regression model \citep{K05}, implemented by the R function {\tt rq} in the package {\tt quantreg}.

{\sf QRCM}: The quantile regression coefficients modeling \citep{fb16} , implemented by the R function {\tt iqr} in the package {\tt qrcm}.

{\sf QPAM} : The quantile partially linear additive model \citep{SL16}, where discrete (continuous) covariates are taken as the linear (nonparametric) parts.

{\sf JQLMA}: The Jackknife quantile linear model averaging \citep{LS15}.

{\sf JCQLMA}: The Jackknife composite quantile linear model averaging \citep{wz19a}.

 It is well known that the performance of model averaging depends on the weight selection criterion, and thus we consider three versions for our proposed procedure.

{\sf EW}: Equal weights are utilized to make predictions.

{\sf QPL}: We randomly set a component of the model weight vector as one and the rest components are taken as zero, indicating that a traditional partially linear model is used for prediction.

{\sf JQPLMA}: The proposed optimal model averaging strategy given in Sub-section \ref{sect.2.2}.

In all simulation examples, we generate a training data set of sample size $n$ to estimate unknown parameters, nonparametric functions and model weights, and generate extra $100$ observations (a testing set) to calculate prediction performances. We use the sample version of ${\rm OAQPE}_n$ (given in Section \ref{sect.3}) to measure accuracy of the out-of-sample prediction performance, defined by
\[{\rm{OAQPE}}=\sum\limits_{k = 1}^n {\sum\limits_{i \in {\cal I}} {\frac{{{\rho _{\tau_k} }\left( {{Y_i} - \hat\mu \left(\bm X_i,\tau _k \right)} \right)}}{{n\left| {\cal I} \right|}}} } ,\]
where $\tau_k=k/(n+1),k=1,\cdots,n$, $\hat\mu \left(\bm X_i,\tau \right)$ is an estimator of the $\tau$th conditional quantile function
$\mu\left(\bm X_i,\tau \right)$ and $\mathcal{I}$ stands for the testing set with the size $\left| {\cal I} \right|$. Following \cite{ls21}, we construct the following three comparison measures
\begin{align}
{\rm Average}~{\rm OAQPE}_{A}=&R^{-1}\sum\limits_{r= 1}^R{\rm{OAQPE}}\left( r \right)_{A}, \nonumber\\
{\rm{Winning ~Ratio}}_{A}=&R^{-1}\sum\limits_{r= 1}^R{\mathbb{I}}\big\{ {\rm{OAQPE}}\left( r \right)_{A}<{\rm{OAQPE}}\left( r \right)_{B},\cdots, \nonumber\\
&{\rm{OAQPE}}\left( r \right)_{A}<{\rm{OAQPE}}\left( r \right)_{H}\big\}, \nonumber\\
\mbox{Loss to}~{\sf JQPLMA}_{A}=&R^{-1}\sum\limits_{r= 1}^R{\mathbb{I}}\left\{ {\rm{OAQPE}}\left( r \right)_{{\sf JQPLMA}}<{\rm{OAQPE}}\left( r \right)_{A}\right\}\nonumber
\end{align}
where $\mathbb{I}\left(\mathcal{C}\right)$ is an indicator function for event $\mathcal{C}$, ${\rm OAQPE}(r)$ is the value of ${\rm OAQPE}$ in the $r$th replication for $r=1,\cdots,R$ and each subscript denotes generic notation for a prediction approach. Please note that
the loss to {\sf JQPLMA} ratio gives us more direct binary comparison of each approach to
{\sf JQPLMA}. Obviously, the smaller ${\rm OAQPE}$ and the bigger ${\rm{winning~ratio}}$, the method is better. Here the total number of replication is taken as $R=200$.

\textbf{Example 1}. In this example, we generate the random samples from the following partially linear additive model
\begin{align} \label{s.1}
{Y_i} = &6X_{i1}+4m_1(X_{i2})+4m_2(X_{i3})+3m_3(X_{i4})+2X_{i5}+2 X_{i6}+2 X_{i7}\nonumber\\
&-2 X_{i8}-2 X_{i9}-2 X_{i10}+\sigma \varepsilon _i , i=1,\cdots,n,
\end{align}
where $m_1\left( u\right)=(2u-1)^2$, $m_2\left( u\right) ={\rm sin}\left(2\pi u\right)/(2-{\rm sin}\left(2\pi u\right))$, $m_3\left( u\right)=0.1{\rm sin}\left(2\pi u\right)+0.2{\rm cos}\left(2\pi u\right)+0.3\left({\rm sin}\left(2\pi u\right)\right)^2+0.4\left({\rm cos}\left(2\pi u\right)\right)^3+0.5\left({\rm sin}\left(2\pi u\right)\right)^3$. The covariates $\bm X_{i}=\left(X_{i1},\cdots,X_{i6}\right)^\top$ are simulated according to ${X_{il}} = \left( {{W_{il}} + t{U_{i}}} \right)/\left( {1 + t} \right)$ for $1\leq l \leq 6$, $X_{i7},X_{i8}\mathop \sim \limits^{i.i.d.} Binomial(1,0.5)$ and $X_{i9},X_{i10}\mathop \sim \limits^{i.i.d.} Binomial(2,0.5)$, where $W_{il}$ and ${U_{i}}$ are generated independently from $Uniform(0,1)$. We have $Corr( X_{il},X_{il'} )=t^2/(1+t^2)$ for $l \neq l'$ by a simple calculation, and set $t=0$, $1$ and 3, representing uncorrelated ($\rho_x=0$), moderate ($\rho_x=0.5$) and high ($\rho_x=0.9$) correlations between covariates. The random error $\varepsilon_i$ is distributed as $N(0,1)$. As in \cite{LS15} and \cite{zw19}, we change the value of $\sigma$ so that the population ${\mathbb{R}^2} = \left\{ {Var\left( {{Y_i}} \right) - Var\left( {\sigma {\varepsilon _i}} \right)} \right\}/Var\left( Y_i \right) = 0.2,0.4,0.6,0.8$ , where $Var\left( \cdot\right)$ represents the sample variance. In this example, only first six covariates are continuous and might be served as the nonparametric component, resulting in $6$ partially linear sub-models for our model averaging procedure. It's worth noting that our mission is to achieve the goal of accurately predicting the joint conditional quantile function $\mu\left(\bm X_i,\tau\right)=6X_{i1}+4m_1(X_{i2})+4m_2(X_{i3})+3m_3(X_{i4})+2X_{i5}+2 X_{i6}+2 X_{i7}-2 X_{i8}-2 X_{i9}-2 X_{i10}+ \sigma Q\left(\varepsilon _i,\tau\right)$ rather than estimate the parameters and nonparametric function in (\ref{s.1}), where $Q\left(\varepsilon _i,\tau\right)$ the $\tau$th quantile function of $\varepsilon _i$.

\textbf{Example 2}. To reflect the flexibility of our procedure, we consider the following multivariate nonparametric regression model with heteroscedasticity
\begin{align}
{Y_i} =&4{\rm{cos}}\left(X_{i1}X_{i2}X_{i3}X_{i4}\right)X_{i5}X_{i6}-3{\rm sin} \left(X_{i7}X_{i8}X_{i9}X_{i10}/4\right)\nonumber\\
&+\left( {\left| {0.5{X_{i9}} - 0.5{X_{i10}}} \right| + 1} \right){\varepsilon _i}, i=1,\cdots,n,
\nonumber
\end{align}
where $\bm X_i=\left(X_{i1},\cdots,X_{ip} \right)^\top$ are generated from a multivariate normal distribution with mean zero and $Cor\left( {{X_{ij}},{X_{il}}} \right) = {0.5^{\left| {j - l} \right|}}$ for $1\leq j, l \leq p$. We also consider $p = 10$ and $15$, corresponding to different sparsity levels.
To assess the robustness and flexibility, we consider six distributions for the random error ${\varepsilon _i}$, including standard normal distribution (\emph{case1}), $t$-distribution with three degrees of freedom (\emph{case2}), a mixture of two normal distributions (\emph{case3}), which is a mixture of $N(0,1)$ and $N(0,25)$ with the weights 95\% and 5\%, $\chi^2$-distribution with one degree of freedom (\emph{case4}), Gamma-distribution $G(1,1)$ (\emph{case5}) and Log normal distribution (\emph{case6}) with the mean and standard deviation of the distribution on the log scale being 0.5 and 0.5, respectively.
In this example, all covariates are continuous and thus any covariate can be taken as the nonparametric component, resulting in $p$ partially linear sub-models for our model averaging procedure.

\tabcolsep=4pt
\begin{table}\scriptsize
\caption{{Simulation results of ${\rm{OAQPE}}$ and $\hat{\bm w}$ for {\sf JQPLMA} with $\mathbb{R}^2=0.8$ and $n=200$ for different correlations and basis functions in example 1.}}
\label{table1}
\begin{tabular}{cccccccccccc} \noalign{\smallskip}\hline
\multicolumn{1}{c}{correlation}
&\multicolumn{1}{c}{basis}
&\multicolumn{1}{c}{${\rm{OAQPE}}$}
&\multicolumn{1}{c}{$\hat w_1$}
&\multicolumn{1}{c}{$\hat w_2$}
&\multicolumn{1}{c}{$\hat w_3$}
&\multicolumn{1}{c}{$\hat w_4$}
&\multicolumn{1}{c}{$\hat w_5$}
&\multicolumn{1}{c}{$\hat w_6$}
\\
\hline
$\rho_x=0$&$\bm b\left( \tau\right)^{(1)}$&0.769(0.058)&0(0)&0.366(0.148)&0.265(0.166)&0.369(0.157)&0(0)&0(0)\\
&$\bm b\left( \tau\right)^{(2)}$&0.775(0.058)&0(0)&0.365(0.143)&0.266(0.160)& 0.368(0.149)&0(0)&0(0.001)\\
&$\bm b\left( \tau\right)^{(3)}$&0.775(0.058)&0(0)&0.365(0.143)&0.267(0.159)& 0.368(0.148)&0(0)&0(0.001)\\
$\rho_x=0.5$&$\bm b\left( \tau\right)^{(1)}$&0.707(0.053)&0(0)&0.249(0.164)&0.174(0.162)&0.574(0.126)&0(0)&0.002(0.024)\\
&$\bm b\left( \tau\right)^{(2)}$&0.712(0.053)&0(0.001)&0.257(0.163)&0.178(0.159)&0.561(0.126)&0(0)&0.003(0.030)\\
&$\bm b\left( \tau\right)^{(3)}$&0.711(0.053)&0(0)&0.259(0.157)&0.181(0.157)&0.557(0.121)&0(0)&0.003(0.028)\\
$\rho_x=0.9$&$\bm b\left( \tau\right)^{(1)}$&0.710(0.057)&0.025(0.058)&0.331(0.138)&0.205(0.138)&0.396(0.124)&0.017(0.045)&0.026(0.063)\\
&$\bm b\left( \tau\right)^{(2)}$&0.715(0.057)&0.028(0.062)&0.330(0.134)&0.206(0.134)&0.385(0.124)&0.021(0.046)&0.030(0.066)\\
&$\bm b\left( \tau\right)^{(3)}$&0.714(0.057)&0.032(0.063)&0.328(0.132)&0.208(0.134)&0.382(0.117)&0.020(0.047)&0.030(0.065)\\
\hline
\end{tabular}
\\
Note: The standard errors of ${\rm{OAQPE}}$ and $\hat{\bm w}=\left(\hat w_1,\cdots,\hat w_6\right)\top$ are denoted inside the parentheses.
\end{table}

\tabcolsep=7pt
\begin{table}\scriptsize
\caption{\leftline{Simulation results over various $\mathbb{R}^2$ and $\rho_x$ in example 1.}}
\label{table2}
\begin{tabular}{cccccccccccc} \noalign{\smallskip}\hline
\multicolumn{1}{c}{correlation}
&\multicolumn{1}{c}{$\mathbb{R}^2$}
&\multicolumn{1}{c}{${\sf QLRM}$}
&\multicolumn{1}{c}{${\sf QRCM}$}
&\multicolumn{1}{c}{${\sf QPAM}$}
&\multicolumn{1}{c}{${\sf JQLMA}$}
&\multicolumn{1}{c}{${\sf JCQLMA}$}
&\multicolumn{1}{c}{${\sf EW}$}
&\multicolumn{1}{c}{${\sf QPL}$}
&\multicolumn{1}{c}{${\sf JQPLMA}$}
\\
\hline
&\multicolumn{9}{c}{\underline{Average ${\rm{OAQPE}}$}}
\\
$\rho_x=0$&0.2&2.544& 2.522& 2.921& 2.522& 2.465& 2.482 &2.515& 2.477\\
&&(0.187)& (0.187)& (0.216)& (0.179)& (0.251)& (0.254)& (0.257)& (0.252)\\
&0.4&1.644& 1.630& 1.938& 1.652& 1.606& 1.594& 1.629& 1.576\\
&&(0.123) &(0.123)& (0.140)& (0.124)& (0.122) &(0.123)& (0.130)& (0.122)\\
&0.6&1.211 &1.201& 1.488& 1.212& 1.181& 1.153& 1.190& 1.117\\
&&(0.089)& (0.089) &(0.110)& (0.087)& (0.089)& (0.085)& (0.089)& (0.083)\\
&0.8&0.901& 0.894 &1.205 &0.904 &0.881& 0.828 &0.871 &0.769\\
&&(0.073) &(0.073) &(0.095) &(0.072)& (0.073)& (0.067)& (0.082)& (0.061)\\
&\multicolumn{9}{c}{\underline{Winning Ratio}}
\\
&0.2&0.0\% &  0.0\% &  0.0\% &  1.0\% & 49.5\% & 13.0\% &  5.0\% & 28.5\% \\
&0.4&0.0\% &  0.0\% &  0.0\% &  1.0\% & 17.5\% & 13.5\% &  6.5\% & 61.5\% \\
&0.6&0.0\% &  0.0\% &  0.0\% &  0.0\% &  2.5\% &  6.0\% &  4.5\% & 87.0\% \\
&0.8&0.0\% & 0.0\% &  0.0\% &  0.0\% &  0.0\% &  0.5\% &  3.0\% & 96.5\% \\
&\multicolumn{9}{c}{\underline{Loss to ${\sf JQPLMA}$}}
\\
&0.2& 88.5\% &77.0\%& 100.0\%& 69.5\% &41.0\% &54.0\% &77.0\%&NA\\
&0.4&96.5\% & 91.5\% & 100.0\% & 91.0\% & 74.5\% & 74.5\% & 86.0\% &NA\\
&0.6&99.5\% &  98.5\% &  100.0\% &  99.5\% &  94.5\% &  91.5\% &  94.5\% &NA\\
&0.8&100.0\% & 100.0\% & 100.0\% & 100.0\% & 100.0\% &  99.5\% &  97.0\%&NA \\
\hline
&\multicolumn{9}{c}{\underline{Average ${\rm{OAQPE}}$}}
\\
$\rho_x=0.5$&0.2&2.528& 2.507 &2.874& 2.513& 2.462 &2.486& 2.520& 2.487\\
&&(0.200)& (0.202)& (0.219)& (0.184)& (0.185)& (0.203)& (0.207)& (0.205)\\
&0.4&1.604 &1.591 &1.878& 1.606 &1.563 &1.564 &1.591& 1.551&\\
&& (0.117) & (0.118) & (0.143) & (0.115) & (0.116) & (0.116) & (0.118) & (0.118)\\
&0.6&1.156& 1.146 &1.394 &1.161& 1.127& 1.110& 1.141& 1.085\\
&& (0.099) & (0.099) &(0.097) & (0.100) & (0.099) & (0.097) & (0.102) & (0.094) \\
&0.8&0.798& 0.792 &1.074& 0.802& 0.780& 0.742& 0.776 &0.707\\
&&(0.063) &(0.063) &(0.083)&(0.063) &(0.062) &(0.058) &(0.062) &(0.056)\\
&\multicolumn{9}{c}{\underline{Winning Ratio}}
\\
&0.2&0.0\% & 0.0\% &  0.0\% &  2.0\% & 58.5\% & 11.5\% &  6.5\% & 18.5\% \\
&0.4&0.0\% &  0.0\% &  0.0\% &  1.0\% & 36.5\% & 13.0\% &  6.0\% & 43.0\% \\
&0.6&0.0\% &  0.0\% &  0.0\% &  0.0\% & 12.1\% & 12.1\% &  4.0\% & 71.7\% \\
&0.8& 0.0\% &  0.0\% &  0.0\% &  0.0\% &  0.5 \% & 1.0\% &  3.5\% & 95.0\% \\
&\multicolumn{9}{c}{\underline{Loss to ${\sf JQPLMA}$}}
\\
&0.2&82.0\% & 69.5\% & 100.0\% & 62.5\% & 34.0\% & 47.0\% & 76.0\%&NA \\
&0.4&90.5\% & 85.5\% & 100.0\% & 86.0\% & 57.0\% & 65.0\% & 82.0\% &NA\\
&0.6&100.0\% &  98.0\% & 100.0\% &  98.0\% &  82.8\% & 82.8\% &  93.9\% &NA\\
&0.8&100.0\% & 100.0\% & 100.0\% &  99.5\% &  98.5\% &  97.0\% &  96.5\% &NA\\
\hline
&\multicolumn{9}{c}{\underline{Average ${\rm{OAQPE}}$}}
\\
$\rho_x=0.9$&0.2&2.570 &2.548& 2.895& 2.546 &2.498 &2.485 &2.529 &2.491\\
&&(0.202) &(0.201) &(0.220) &(0.194) &(0.194) &(0.196) &(0.200) &(0.194)\\
&0.4& 1.646& 1.632& 1.920 &1.646& 1.603& 1.539 &1.588& 1.542\\
&&(0.128) &(0.129) &(0.155) &(0.129) &(0.126) &(0.124) &(0.132) &(0.128)\\
&0.6&1.227 &1.217& 1.445& 1.234& 1.199& 1.086 &1.137 &1.080\\
&&(0.091) &(0.092) &(0.117) &(0.090) &(0.090) &(0.085) &(0.094) &(0.084)\\
&0.8&0.917& 0.910 &1.139& 0.922& 0.899& 0.726 &0.799 &0.710\\
&&(0.064) &(0.064) &(0.102) &(0.065) &(0.065) &(0.055) &(0.065) &(0.057)\\
&\multicolumn{9}{c}{\underline{Winning Ratio}}
\\
&0.2&0.0\% &  0.0\% &  0.0 \% & 0.5 \% &35.0 \% &30.0 \% & 8.5 \% &23.5\% \\
&0.4&0.0\% &  0.0\% &  0.0\% &  0.0\% & 10.5\% & 39.0\% & 11.5\% & 39.0\% \\
&0.6&0.0\% &  0.0\% &  0.0\% &  0.0\% &  0.5\% & 39.0\% &  4.0\% & 56.5\% \\
&0.8&0.0\% &  0.0\% &  0.0\% &  0.0\% &  0.0\% & 25.5\% &  1.5\% & 73.0\% \\
&\multicolumn{9}{c}{\underline{Loss to ${\sf JQPLMA}$}}
\\
&0.2&89.5\% & 81.5\% & 100.0\% & 77.0\% & 50.5 \% &37.5\% & 76.0\% &NA\\
&0.4&93.0\% & 91.5\% & 100.0\% & 90.5\% & 83.5\% & 44.0\% & 81.5\% &NA\\
&0.6&100.0\% &  99.0\% &  100.0\% &  99.5\% &  98.5\% &  59.5 \% & 92.0\% &NA\\
&0.8&100.0\% & 100.0\% & 100.0\% & 100.0\% & 100.0\% &  74.0\% &  98.5\% &NA\\
\hline
\end{tabular}
\\
Note: The standard error of ${\rm{OAQPE}}$ is denoted inside the parentheses.
\end{table}

\tabcolsep=7pt
\begin{table}\scriptsize
\caption{\leftline{Simulation results over various errors and sparsity levels in example 2.}}
\label{table3}
\begin{tabular}{cccccccccccc} \noalign{\smallskip}\hline
&\multicolumn{1}{c}{error}
&\multicolumn{1}{c}{${\sf QLRM}$}
&\multicolumn{1}{c}{${\sf QRCM}$}
&\multicolumn{1}{c}{${\sf QPAM}$}
&\multicolumn{1}{c}{${\sf JQLMA}$}
&\multicolumn{1}{c}{${\sf JCQLMA}$}
&\multicolumn{1}{c}{${\sf EW}$}
&\multicolumn{1}{c}{${\sf QPL}$}
&\multicolumn{1}{c}{${\sf JQPLMA}$}
\\
\hline
&\multicolumn{9}{c}{\underline{Average ${\rm{OAQPE}}$}}
\\
$p=10$&\emph{case 1}&1.170& 1.162& 1.134& 1.132& 1.119& 1.095 &1.157& 0.964\\
&&(0.160) &(0.161)& (0.152)& (0.158)& (0.158)& (0.149)& (0.168)& (0.119)\\
&\emph{case 2}&1.298 &1.287 &1.271& 1.255& 1.238& 1.221& 1.279& 1.100\\
&&(0.164)& (0.164)& (0.155)& (0.159) &(0.158)& (0.156)& (0.173)& (0.138)\\
&\emph{case 3}&1.241& 1.234 &1.186 &1.198& 1.183& 1.167 &1.230 &1.049\\
&&(0.144)& (0.146)& (0.148)& (0.145)& (0.145)& (0.138)& (0.152)& (0.120)\\
&\emph{case 4}&1.222& 1.212 &1.255& 1.177 &1.161& 1.145& 1.195 &1.015\\
&&(0.136)& (0.138)& (0.153)& (0.135)& (0.136)& (0.130)& (0.152)& (0.114)\\
&\emph{case 5}&1.159& 1.151& 1.221& 1.121& 1.107& 1.083 &1.151 &0.944\\
&&(0.144)& (0.145)& (0.159)& (0.142)& (0.141)& (0.134)& (0.152)& (0.113)\\
&\emph{case 6}&1.177& 1.169& 1.386& 1.136& 1.121 &1.101& 1.159 &0.971\\
&&(0.146)& (0.145)& (0.178)& (0.141)& (0.140)& (0.137)& (0.157)& (0.113)\\
&\multicolumn{9}{c}{\underline{Winning Ratio}}\\
&\emph{case 1}&0.0\% & 0.0\% & 0.0\%&  0.0\%&  3.5\%&  1.0\%&  0.5\% &95.0\%\\
&\emph{case 2}&0.0\% &0.0\%  &0.0\% & 0.0\% &4.0\%&  2.5\% & 0.5\% &93.0\%\\
&\emph{case 3}&0.0\% &  0.0\% &  0.0\% &  0.0\% &  5.0\% &  1.5\% &  1.0\% & 92.5\% \\
&\emph{case 4}&0.0\% &  0.0\% &  0.0\% &  0.0\% &  2.0\% &  1.0\% &  0.5\% & 96.5\% \\
&\emph{case 5}&0.0\% &  0.0\% &  0.0\% &  0.0\% &  2.0\% &  0.5\% &  0.5\% & 97.0\% \\
&\emph{case 6}&0.0\% &  0.0\% &  0.0\% &  0.0\% &  3.0\% &  1.5\% &  1.5\% & 94.0\% \\
&\multicolumn{9}{c}{\underline{Loss to ${\sf JQPLMA}$}}\\
&\emph{case 1}&99.5\%& 99.5\%& 100.0\% &97.0\%& 96.0\% &97.0\%& 98.0\%&NA\\
&\emph{case 2}&99.5\%& 99.5\%& 100.0\%& 97.0\%& 94.0\%& 95.0\% &99.0\%&NA\\
&\emph{case 3}&98.5\% & 97.5\% & 100.0\% & 95.5 \% &94.0\% & 95.0\% & 98.5\% &NA\\
&\emph{case 4}& 99.0\% & 98.5\% & 100.0\% & 98.5\% & 97.5\% & 97.0\% & 99.0\% &NA\\
&\emph{case 5}& 98.0\% & 98.0\% & 100.0\% & 98.0\% & 97.5\% & 98.0\% & 98.5\% &NA\\
&\emph{case 6}&99.0\% & 98.5\% & 100.0\% & 97.0\% & 96.0\% & 96.0\% & 97.0\% &NA\\
\hline
&\multicolumn{9}{c}{\underline{Average ${\rm{OAQPE}}$}}
\\
$p=15$&\emph{case 1}&1.213& 1.209& 1.237& 1.143 &1.125& 1.153& 1.203& 1.000\\
&&(0.147)& (0.151)& (0.162)& (0.148)& (0.147)& (0.145)& (0.152)& (0.126)\\
&\emph{case 2}&1.343 &1.339& 1.411& 1.267 &1.243 &1.284 &1.353 &1.140\\
&&(0.156)& (0.163)& (0.172)& (0.154)& (0.155)& (0.156)&(0.167)& (0.143)\\
&\emph{case 3}&1.276& 1.271& 1.315 &1.205 &1.184& 1.217 &1.281& 1.062\\
&&(0.151)& (0.156)& (0.173)& (0.148)& (0.150)& (0.149)& (0.162)& (0.129)\\
&\emph{case 4}&1.278 &1.276& 1.388& 1.208& 1.186& 1.221& 1.288& 1.069\\
&&(0.151)& (0.155)& (0.173)& (0.144)& (0.142)& (0.147)& (0.155)& (0.124)\\
&\emph{case 5}&1.182& 1.182& 1.359& 1.116 &1.095& 1.127& 1.190 &0.973\\
&& (0.144)& (0.147)& (0.177)& (0.140)& (0.141)& (0.139)& (0.153)& (0.113)\\
&\emph{case 6}&1.215& 1.211 &1.465& 1.149& 1.129& 1.156 &1.222 &1.001\\
&&(0.137)& (0.138)& (0.197)& (0.133)& (0.134)& (0.131)& (0.150)& (0.110)\\
&\multicolumn{9}{c}{\underline{Winning Ratio}}\\
&\emph{case 1}&0.0\% &  0.0\% &  0.0\% &  0.5\% &  4.0 \% & 0.5\% &  0.0\% & 95.0\% \\
&\emph{case 2}&0.0\% &  0.0\% & 0.0\% &  0.0\% &  9.0\% &  0.0\% &  0.0\% &  91.0\% \\
&\emph{case 3}&0.0\% &  0.0\% &  0.0\% &  0.0\% &  6.5\% &  0.0\% &  0.0\% & 93.5\% \\
&\emph{case 4}&0.0\% &  0.0\% &  0.0 \% & 0.0\% &  7.0\% &  0.0\% &  0.5\% & 92.5\% \\
&\emph{case 5}&0.0\% &  0.0\% &  0.0\% &  0.0\% &  7.5\% &  0.0\% &  0.5\% & 92.0\% \\
&\emph{case 6}&0.0\% &  0.0\% &  0.0\% &  0.0\% &  8.0\% &  0.0\% &  0.0\% & 92.0\%
\\
&\multicolumn{9}{c}{\underline{Loss to ${\sf JQPLMA}$}}\\
&\emph{case 1}&99.5\% & 99.0\% & 100.0\% & 97.0\% & 95.5\% & 98.5\% & 99.5\% &NA\\
&\emph{case 2}&100.0\% & 100.0\% & 100.0 \% & 97.5\% &  91.0\% & 100.0 \% &100.0\% &NA\\
&\emph{case 3}&100.0\% & 100.0\% & 100.0 \% & 98.0\% &  93.5\% &  99.5\% & 100.0\% &NA\\
&\emph{case 4}& 99.5\% & 99.5\% & 100.0\% & 95.5\% & 93.0 \% &98.0\% & 99.5\% &NA\\
&\emph{case 5}&99.5\% & 98.5\% & 100.0\% & 94.0\% & 92.5\% & 97.0 \% &99.0\% &NA\\
&\emph{case 6}&98.5\% & 99.0\% & 100.0\% & 94.0\% & 92.0\% & 96.5 \% &99.5\% &NA\\
\hline
\end{tabular}
\\
Note: The standard error of ${\rm{OAQPE}}$ is denoted inside the parentheses.
\end{table}

To implement our procedure, we need to determine the degree of B-spline and the number of knots, which play important roles in numerical studies. Recent researching findings (\cite{hwz04} and \cite{k07}) have showed that lower order splines might be better choice, such as linear splines ($d = 2$). It is well-known that higher order splines would induce complicated interactions and collinearity among the variables in the model as the effect of the splines on the model is multiplicative. Therefore, we suggest using linear splines in our simulations because of its desirable properties such as optimality \cite{knp94}. Moreover, we set the number of interior knots as ${N_n} = \left[ n^{1/5} \right]$ with $[s]$ being the largest integer not greater than $s$. In addition, it is a natural question whether our proposed method is sensitive to the choice of basis set $\bm b\left( \tau\right)$. So we conduct a sensitivity analysis for the choice of $\bm b\left( \tau\right)$. Similar to \cite{fb16} and \cite{YCC17}, we consider the following three types of basis functions
\begin{align}
\bm b\left( \tau\right)^{(1)}&=\left(1,\Phi^{-1}\left(\tau\right)\right)^\top,\nonumber\\
\bm b\left( \tau\right)^{(2)}&=\left(1,\tau,\tau^2,\tau^3\right)^\top,\nonumber\\
\bm b\left(\tau\right)^{(3)}&= \left(1,\tau,\Phi^{-1}\left(\tau\right),-\mbox{log}\left(1-\tau\right) \right)^\top,\nonumber
\end{align}
where $\Phi\left(\cdot\right)$ denotes the distribution function of the standard normal distribution. Table \ref{table1} lists the average of ${\rm{OAQPE}}$ and estimated model weight vector $\hat{\bm w}$ for {\sf JQPLMA} with different basis functions for $\mathbb{R}^2=0.8$ and $n=200$ in example 1.
Table \ref{table1} shows the second, third and fourth sub-models carry almost all the weights, and the combination of the three models is indeed the true model, which indicates the proposed cross-validation based
method works very well for selection of the weights in the model averaging prediction.
Furthermore, it is easy to see that there is little difference for ${\rm{OAQPE}}$ and $\hat{\bm w}$ among different basis functions, indicating that our proposal is not sensitive to the selection of $\bm b\left( \tau\right)$. The results of ${\rm{OAQPE}}$ and $\hat{\bm w}$ for other settings are also not sensitive to the selection of $\bm b\left( \tau\right)$. To save space, we don't report the results for other settings. Therefore, we fix $\bm b\left( \tau\right) =\left(1,\Phi^{-1}\left(\tau\right)\right)^\top$ in the simulation studies.

The simulation results over all designs are reported in Tables \ref{table2}--\ref{table3}. We might obtain the following conclusions. Firstly, we check the performance over different signals ($\mathbb{R}^2$) and different levels of dependency among covariates ($\rho_x$). From Table \ref{table2}, we confirm that our proposed approach {\sf JQPLMA} yields the smallest ${\rm{OAQPE}}$ and the highest winning ratio when $\mathbb{R}^2$ varies form 0.4 to 0.8. When $\mathbb{R}^2=0.2$, {\sf JCQLMA} and {\sf EW} are slightly better than {\sf JQPLMA}. One possible explanation is that our procedure requires the estimation of more parameters and nonparametric functions and thus might result in poor estimators for the relatively small signal $\mathbb{R}^2$. It is also interesting that the winning ratio of our method and the loss to {\sf JQPLMA} increase quickly as $\mathbb{R}^2$ increases, indicating that the superiority of {\sf JQPLMA} is increasingly apparent for the large $\mathbb{R}^2$. Secondly, we study the performance over different sparsity levels and error distributions. Table \ref{table3} reveals that {\sf JQPLMA} outperforms all competing methods uniformly. It may not be surprising to understand that traditional model-based approaches ({\sf QLRM}, {\sf QRCM}, {\sf QPAM} and {\sf QPL}) have poor prediction performance because they adopt a single misspecified model structure to make predictions. Furthermore, although all model averaging approaches ({\sf JQLMA}, {\sf JCQLMA}, {\sf EW} and {\sf JQPLMA}) employ misspecified candidate models, the proposed {\sf JQPLMA} optimally combines useful
information from more flexible semiparametric sub-models, and thus produces more accurate prediction performance.

In sum, simulation studies show that our proposed procedure has satisfactory finite sample properties for various settings.

\bibliography{ref}
\bibliographystyle{apalike}

\end{document}